\documentclass[12pt,reqno]{amsart}

\usepackage[dvips]{graphicx} 

\usepackage{
hyperref
}

\usepackage{breakurl}   

\usepackage{color}
\definecolor{darkblue}{rgb}{0.0, 0.0, 0.8}

\newtheorem{theorem}{Theorem}[section]

\theoremstyle{definition}

\numberwithin{equation}{section}

\newcommand\R {{\mathbb R}}
\newcommand\C {{\mathbb C}}

\newcommand\CP {{\mathbb C\mathbb P}} 
\newcommand\RP {{\mathbb R\mathbb P}}

\title 
[Mutually unbiased bases via complex projective trig]
{Mutually unbiased bases via complex projective trigonometry}

\author[Mikhail Katz]{Mikhail G. Katz}\address{M.~Katz, Department of
  Mathematics, Bar Ilan University, Ramat Gan 5290002
  Israel}\email{katzmik@math.biu.ac.il}

\date{\today}


\begin{document}

\thispagestyle{empty}


\begin{abstract}
We give a synthetic construction of a complete system of mutually
unbiased bases in~$\C^3$.
\end{abstract}


\maketitle

\section{Synthetic description}
\label{s1}

In quantum mechanics, the concept of mutually unbiased bases
\cite{Be17} is of fundamental importance.  Recall that orthonormal
bases~$(e_i)$ and~$(f_j)$ for~$\C^n$ are called \emph{mutually
unbiased} if the absolute value of the Hermitian inner product
$\langle e_i, f_j\rangle$ for all~$i,j$ is independent of~$i$ and~$j$
(and therefore necessarily equals~$\frac{1}{\sqrt{n}}$).  No knowledge
of quantum theory is required for reading this article.

{A \emph{complete system} of mutually unbiased bases in $\C^n$ is
  a system of $n+1$ such bases}.  It is still unknown whether a
complete system of bases exists for~$\mathbb C^6$ i.e., in~$\mathbb
C\mathbb P^5$.  A complete system of mutually unbiased bases is known
to exist when the vector space dimension is a prime power.  We propose
a new approach to constructing mutually unbiased bases, and implement
it in the case of the complex projective plane.

Fix a complex projective line
\begin{equation}
\label{e11}
\ell\subseteq \CP^2.
\end{equation}
{Recall that $\ell$ is isometric to a 2-sphere of constant
  Gaussian curvature.}  Let
\[
A,B,C,D\in\ell
\]
be the vertices of a regular inscribed tetrahedron, {i.e., an
  equidistant 4-tuple of points on the 2-sphere.}  There is a
synthetic construction of a complete system of mutually unbiased bases
(MUBs), where each basis includes one of the 4 points $A,B,C,D$.  Let
\begin{equation}
\label{e12}
m\subseteq\CP^2
\end{equation}
be the complex projective line consisting of points at maximal
distance from~$A\in\CP^2$.  On the projective line $\ell$ of
\eqref{e11}, let
\[
E\in\ell
\]
be the antipodal point of~$A$, so that~$E=\ell\cap m$, where $m$ is
the complex projective line of \eqref{e12}.  {Then $\{A,E\}$ is a
  pair of points at maximal geodesic distance in $\CP^2$.  In
  homogeneous coordinates, this corresponds to the fact that their
  representing vectors in $\C^3$ are orthonormal with respect to the
  standard Hermitian inner product.}

Choose a great circle (i.e., closed geodesic)
\[
S^1\subseteq m
\]
passing through~$E$.  We view~$S^1$ as the equator of~$m$ and denote
by
\[
A',A''\in m
\]
{the corresponding north and south poles}.  Then~$S^1\subseteq m$
is the set of equidistant points from~$A'$ and~$A''$.  Let
\[
\gamma_{ABE}
\]
be the closed geodesic (great circle) passing through the three points
$A,B,E\in\ell$.  Consider the real projective
plane
\[
\RP^2_B\subseteq\CP^2
\]
that includes the great circles~$\gamma_{ABE}$ and~$S^1$.  Consider
also the real projective plane~$\RP^2_C\subseteq\CP^2$ that includes
the great circles~$\gamma_{ACE}$ and~$S^1$.

\begin{theorem}
\label{t11}
A complete system of MUBs is constructed as follows.  The first basis
is~$\{A,A',A''\}\subseteq\CP^2$.  {On the complex projective line
  passing through $A'$ and $A''$, let $S^1$ be the equidistant circle
  from $A'$ and $A''$.  Choose an equilateral triangle~$EE'E''$
  in~$S^1$.  Let $\ell$ be the projective line through $A$ and $E$.
  We complete $A$ to an equidistant $4$-tuple $A,B,C,D$ on $\ell$.}
Next, we complete~$B$ to a basis
$\{B,B',B''\}\subseteq\RP^2_B\subseteq\CP^2$ by choosing~$B'$
appropriately on the great circle~$\gamma_{AE'}\subseteq\RP^2_B$, and
similarly for~$B''\in \gamma_{AE''}$.  Next, we complete~$C$ to a
basis~$\{C,C',C''\}$ as we did for~$B$; the same for~$D$.
\end{theorem}
The remainder of the paper presents a proof of Theorem~\ref{t11}.

\section{Complex projective trigonometry}

Consider a pair of unit vectors~$v,w$ in~$\C^n$ endowed with its
Hermitian inner product~$\langle\;,\;\rangle$.  Associated with the
pair~$v,w$, there is a pair of angles:
\begin{enumerate}
\item
$\alpha= \arccos \text{Re} \langle v,w\rangle$ is the usual angle
  between the vectors when~$\C^n$ is identified with~$\R^{2n}$, and
\item
$\theta= \arccos |\langle v,w\rangle|$ is the least angle between
  vectors in the complex lines spanned by~$v$ and~$w$.
\end{enumerate}
Letting~$Pw$ be the orthogonal projection of~$w$ to the complex line
spanned by~$v$, we denote by
\begin{equation}
\label{e21b}
\psi
\end{equation}
the angle between~$v$ and~$Pw$, so that by the spherical law of
cosines {(see \cite[p.\,17]{Br99})} we have
\begin{equation}
\label{e21}
\cos\alpha=\cos\theta\cos\psi.
\end{equation}

Given a pair of geodesics issuing from a point~$A\in\CP^n$, the
corresponding angle~$\alpha$ (resp.~$\theta$) is defined similarly
using their unit tangent vectors at~$A$.  Consider a geodesic triangle
in~$\CP^n$ with sides of length~$a,b,c$ and angles~$\alpha$ and
$\theta$ at the point opposite the side~$c$.  We normalize the metric
on~$\CP^n$ so that the sectional curvature satisfies~$\frac14\leq
K\leq1$.  Then each complex projective line~$\CP^1\subseteq\CP^n$ is
isometric to a unit sphere.  The following formula, generalizing the
law of cosines of spherical trigonometry, goes back to Shirokov
\cite{Sh57} and was exploited in \cite{Ka84} in 1984 as well as in
1991 in \cite[p.\,176]{Ka91}:
\begin{equation}
\label{e22}
\cos c = \cos a \cos b + \sin a \sin b
\cos\alpha-2\sin^2\frac{a}2\sin^2 \frac{b}2 \sin^2\theta.
\end{equation}
Let~$A\in\CP^2$ and let~$m\subseteq\CP^2$ be the complex projective
line consisting of points at maximal distance~$\pi$ from~$A$ (the
notation is consistent with that introduced in Section~\ref{s1}).

Let~$E,X\in m$ and consider a geodesic triangle~$AEX$ contained in a
copy of a real projective plane (of curvature $\frac14$ according our
normalisation).  As before,~$\alpha$ and~$\theta$ are the angles
between the tangent vectors at~$A$ to the chosen
geodesics~$\gamma_{AE}$ and~$\gamma_{AX}$.

Consider the tangent vector of the geodesic~$\gamma_{AX}$ at the
point~$X$.  If we parallel translate this vector along the
geodesic~$\gamma_{XE}$ to a vector~$u$ at the point~$E$, then
\begin{enumerate}
\item
the angle at~$E$ between~$u$ and the tangent vector to the
geodesic~$\gamma_{AC}$ is the angle~$\psi$ of \eqref{e21b} by
\cite[p.\,177]{Ka91};
\item
the distance~$d(E,X)$ between~$E$ and~$X$ is~$d(E,X)=2\theta$
\end{enumerate}
(the factor of~$2$ is due to our normalisation of the metric).

We complete~$A$ to a basis~$\{A,A',A''\}$ by choosing a pair of
antipodal points~$A',A''\in m$.

The following identity plays a key role in the construction.  Relative
to the above normalisation of curvature, a totally geodesic~$\CP^1$
has Gaussian curvature~$1$ whereas a totally geodesic~$\RP^2$ has
Gaussian curvature~$\frac14$.  Let~$d$ be the spherical sidelength of
a regular inscribed tetrahedron in~$\CP^1$.  Then geodesic arcs of
length~$d$,~$\pi-d$, and~$\frac\pi2$ form a right-angle triangle in
$\RP^2$ with hypotenuse~$d$, so that
\begin{equation}
\label{e23}
\cos\tfrac{d}2=\cos\big(\tfrac{\pi-d}2\big)\cos\tfrac\pi4.
\end{equation}

\section{Mutually unbiased bases}
\label{s3}

By definition, another basis~$\{B,B',B''\}$ is unbiased with respect
to~$\{A,A',A''\}$ if and only if the distance~$d$ between~$A$ and~$B$
is~$d=2\arccos\frac1{\sqrt3}$ (again the factor of~$2$ is due to our
normalisation), and similarly for all the other distances between a
member of the first basis and a member of the second basis.  Note
that~$\cos d=2\cos^2\frac{d}2-1=\frac23-1=-\frac13$, i.e.,
$d=\arccos\big(\hspace{-2pt}-\frac13\big)$ is the spherical sidelength of a
regular inscribed tetrahedron.

Consider the equator~$S^1\subseteq m$ equidistant from~$A'$ and~$A''$,
and choose an equilateral triangle~$EE'E''$ on~$S^1$, so that
$d(E,E')=\frac{2\pi}3$ and similarly for the other two pairs.  The
complex projective line through~$A$ and~$E$ is denoted~$\ell$ as in
Section~\ref{s1}.

We choose a geodesic arc~$\gamma_{AE}\subseteq\ell\subseteq\CP^2$ and
a point~$B\in \gamma_{AE}$ so that~$d(A,B)=d$.  Next, we consider the
complex projective line through~$A$ and~$E'$.  On this line, we choose
a geodesic~$\gamma_{AE'}\subseteq\CP^2$ so that the angle~$\psi$
between~$\gamma_{AE}$ and~$\gamma_{AE'}$ is~$\psi=\pi$
(here~$\alpha=\frac{2\pi}3$ and~$\theta=\frac\pi3$).  Then the
geodesics~$\gamma_{AE}$ and~$\gamma_{AE'}$ lie in a common real
projective plane~$\RP^2$ which also includes the equator~$S^1\subseteq
m$.  We choose the point~$B'\in \gamma_{AE'}$ so that~$d(A,B')=d$.  We
choose the point~$B''$ above~$E''$ similarly.  Note that the
angle~$\psi$ between~$\gamma_{AE'}$ and~$\gamma_{AE''}$ is also~$\pi$,
because parallel transport around~$S^1$ in~$\RP^2$ of a vector~$u$
orthogonal to~$S^1$ gives~$-u$.

By construction,~$d(A,B)=d$.  Let us calculate~$d(A',B)$.  Consider
the right-angle triangle~$A'EB$ in~$\RP^2$.  As
$\cos\frac{d}2=\frac1{\sqrt{3}}$ by the unbiased condition, we have
$\sin\frac{d}2=\sqrt{\frac23}$.  Hence by the spherical theorem of
cosines,
\[
\cos(\tfrac12 d(A',B))= \cos(\tfrac\pi4)\cos(\tfrac\pi2 -
\tfrac{d}2)= \tfrac1{\sqrt{2}} \sin\tfrac{d}2 = \tfrac1{\sqrt{2}}
\sqrt{\tfrac23}=\tfrac1{\sqrt{3}},
\]
as required.  The equality~$AB=A'B$ is equivalent to \eqref{e23}.  A
similar calculation shows that~$\cos d(A'',B)=-\frac13$, as well.  By
symmetry, each of the points~$B'$ and~$B''$ is also at geodesic
distance~$d$ from~$A$,~$A'$ and~$A''$.  It remains to check that the
triple~$\{B,B',B''\}$ corresponds to an orthonormal basis.
Indeed,~$\psi=\pi$ implies~$\cos\alpha=\cos\psi\cos\theta=-\cos\theta$
by \eqref{e21}.  We apply \eqref{e22} at the point~$A$ we obtain
\[
\begin{aligned}
\cos d(B,B') &= \cos^2 d + \sin^2 d \cos\alpha -
2\sin^4\tfrac{d}2\sin^2\theta \\&=
\tfrac19+\tfrac89(-\tfrac12)-2\cdot\tfrac49\cdot\tfrac34=-1
\end{aligned}
\]
and therefore~$d(B,B')=\pi$ as required.

\section{Completing the MUBs}

To specify a third basis, we choose~$C\in\ell$ so as to complete~$A$
and~$B$ to an equilateral triangle~$ABC$ with side~$d$ (and angle
$\frac{2\pi}3$).  Then the corresponding angle~$\theta$ vanishes, and
hence~$\alpha=\psi=\frac{2\pi}3$.  We complete~$C$ to a basis
$\{C,C',C''\}$ as in Section~\ref{s3}, using a real projective plane
which includes the geodesic~$\gamma_{ACE}$ and the
equator~$S^1\subseteq m$ (the same equator as before, equidistant
between~$A'$ and~$A''$).  To check that the~$C$-basis is unbiased
relative to the~$B$-basis, note that~$\psi(B,B')=\pi$ and
$\psi(B,C)=\frac{2\pi}3$ hence
$\psi(B',C)=\pi-\frac{2\pi}3=\frac\pi3$, while~$\theta=\frac\pi3$ as
before.  Hence~$\cos d(B,C')=\frac19 +\frac89 \cos\psi\cos\theta -
2\cdot\frac49\sin^2\theta= -\frac13$, as required.  The other pairs
are checked similarly.

To specify a fourth basis, we choose~$D\in\ell$ to complete the
equilateral triangle~$ABC$ to a regular inscribed tetrahedron, and
proceed as before.  This completes the proof of Theorem~\ref{t11}.

Such a synthetic construction of MUBs involves less redundancy than
complex Hadamard matrix presentations, and may therefore be useful in
shedding light on MUBs in higher-dimensional complex projective
spaces, such as~$\CP^5$ where the maximal number of MUBs is unknown.

\section{Funding information}

Mikhail Katz was supported by the BSF grant 2020124 and the ISF grant
743/22.

\end{document}